\documentclass[preprint,authoryear,12pt]{article}

\usepackage{amsmath}
\usepackage{amssymb}
\usepackage{multirow}
\usepackage[english]{babel}
\usepackage{graphicx}
\usepackage{color}
\usepackage{natbib}
\usepackage{amssymb}
\usepackage[margin=1 in]{geometry}
\usepackage{setspace}
\usepackage{enumerate}
\renewcommand{\vec}[1]{\mathbf{#1}}
\newcommand{\svec}[1]{\boldsymbol{#1}}

\newcommand{\yv}{\vec{y}}
\newcommand{\yvb}{\bar{\yv}}

\newcommand{\xv}{\vec{x}}
\newcommand{\xiv}{\svec{\xi}}
\newcommand{\zetav}{\svec{\zeta}}
\newcommand{\Zv}{\vec{Z}}
\newcommand{\zv}{\vec{z}}
\newcommand{\cv}{\vec{c}}
\newcommand{\uv}{\vec{u}}
\newcommand{\Cv}{\vec{C}}

\newcommand{\Pv}{\vec{P}}

\newcommand{\Sv}{\vec{S}}
\newcommand{\etav}{\svec{\eta}}
\newcommand{\muv}{\svec{\mu}}

\newcommand{\betah}{\hat{\beta}}
\newcommand{\betav}{\svec{\beta}}
\newcommand{\betavh}{\hat{\svec{\beta}}}
\newcommand{\alphah}{\hat{\alpha}}
\newcommand{\alphav}{\svec{\alpha}}
\newcommand{\deltav}{\svec{\delta}}
\newcommand{\Lambdav}{\svec{\Lambda}}
\newcommand{\Omegav}{\svec{\Omega}}
\newcommand{\hv}{\vec{h}}
\newcommand{\Hv}{\vec{H}}
\newcommand{\gv}{\vec{g}}
\newcommand{\Gv}{\vec{G}}
\newcommand{\piv}{\svec{\pi}}

\newcommand{\av}{\vec{a}}

\newcommand{\bvh}{\hat{\svec{\beta}}}
\newcommand{\etah}{\hat{\etav}}

\newcommand{\xvh}{\hat{\xv}}
\newcommand{\deltavh}{\hat{\deltav}}
\newcommand{\alphavh}{\hat{\alphav}}
\newcommand{\thetav}{\svec{\theta}}

\newcommand{\sv}{\vec{s}}

\newcommand{\Sigmav}{\svec{\Sigma}}
\newcommand{\Sigmavh}{\hat{\Sigmav}}
\newcommand{\fv}{\vec{f}}

\newcommand{\Jv}{\vec{J}}

\newcommand{\Iv}{\vec{I}}

\newcommand{\bv}{\vec{b}}

\numberwithin{equation}{section}

\begin{document}
\doublespacing
\title{On generalized multinomial models and joint percentile estimation}

\author{\small I.~Das$$, S.~Mukhopadhyay$$\thanks{{Corresponding author. Email}:
siuli@math.iitb.ac.in \vspace{6pt} {Phone}: 912225767495}\\\small $^{}$Department of Mathematics, Indian
Institute of Technology Bombay, Mumbai, India}

\date{}
\maketitle
\begin{abstract}
This article proposes a family of link functions for the multinomial response model. The link family
includes the multicategorical logistic link as one of its members. Conditions for the local orthogonality of the
link and the regression parameters are given. It is shown that local orthogonality of the parameters in a
neighbourhood makes the link family location and scale invariant. Confidence regions for jointly
estimating the percentiles based on the parametric family of link functions are also determined. A numerical
example based on a combination drug study is used to illustrate the proposed parametric link family and the confidence
regions for joint percentile estimation.

\end{abstract}

Keywords: confidence regions, multicategorical logistic link, parameter orthogonality, standardization




\section{Introduction}

In this article we address two issues related to multinomial response models, (i) a family of link functions and
(ii) percentile estimation under a parametric link family. In the first few sections we propose  a family of link
functions for multinomial nominal response models. When working with multinomial data sets the common practice is
to fit the multicategory logistic link function \citep[pp.~267-274]{agresti_2002}. However, \cite{1992_czado}
show that if the link function is incorrectly assumed then it leads to biased estimates thus increasing the mean
squared error of prediction. Using a data set based on a combination drug therapy experiment we show that
parameter estimation is improved by using the proposed link family instead of the commonly
used multivariate logistic link. The parametric link family proposed includes the multivariate logistic link as one of its members. %
 In the later part of the article we discuss three methods for finding confidence regions for the percentiles of
a multinomial response model. The confidence regions determined are
based on the estimated values of the link parameters. 

In univariate generalized linear models (GLMs), especially for binary data,  family of link functions have been
discussed by many researchers. Some of the one and two parameter link families for binary models are, proposed by
namely \cite{1975_prentice,1976_prentice,1980_pregibon,1982_johnson,1981_aranda,
1988_stukel,1992_czado,1992_czadob,1993_czado,1997_czado,1999_lang}. However, unlike binary regression models
research papers on link families for multinomial responses are rarely found in the literature. The two parametric
link families proposed by \cite{1985_genter} and \cite{1999_lang} are only applicable to multinomial data sets
with ordered categories. Till date we were unable to find any work which addresses a family of link functions for
multinomial data sets with nominal responses. The situation is similar for percentile estimation methods in
multinomial response models. Though a huge number of research papers (namely,
\cite{1979_hamilton,1986_carter,1986_williams,2001_huang,2006_biedermann,2011_li})
 have been
written on percentile estimation and effect of link misspecification on percentile estimation for binary data,
almost no work has been done in the case of multinomial data. There are, however, many experimental situations
where multinomial responses may be observed for each setting of a group of control variables. As a typical
example we may consider a drug testing experiment, where both the efficacious and toxic responses of the drug/s
are measured on the subjects. This results in two responses, efficacy and toxicity of the drug, both of which are
binary in nature. Since the responses come from the same subject they are assumed to be correlated, and can be
modeled using a multinomial distribution \citep{mukhopadhyaykhuri_2008b}. In this situation it may be of interest to the
experimenter to jointly estimate the $100p$ percentile of the efficacy and toxic responses. In this article we
discuss a numerical example based on the pain relieving and toxic effects of two analgesic drugs and determine confidence
regions for the $100$p percentiles of both the responses.

While parametric link families are able to improve the maximum likelihood fit when compared to canonical links,
any correlation between the link and the regression parameters leads to an increase in the variances of the
parameter estimates [\cite{1997_czado}]. However, it can be shown that if the parameters are orthogonal to each
other then the variance inflation reduces to zero for large sample sizes \citep{1987_cox}. Conditions for local
orthogonality in a neighbourhood was proposed by \cite{1997_czado} for univariate GLMs. In this article we extend
these conditions so that we can
 apply them to a multiresponse situation. It is also shown that the local orthogonality of the
parameters imply location and scale invariance of the family of link functions.

The remainder of the article is organized as follows: In Section \ref{mvglm} we describe the family of
link functions for the multinomial model. Detailed conditions of local orthogonality between the link and the
regression parameters are given in Section \ref{ortho}. In Section \ref{cr} we discuss three interval methods for
percentile estimation in a multinomial model. The proposed link family and confidence regions are illustrated
with a numerical example based on a drug testing experiment in Section \ref{example}. Concluding remarks are
given in Section \ref{conclusion}.

\section{A family of link functions for multinomial data}\label{mvglm}
In this section the multinomial response model with a parametric link function is introduced.
We use a scaled version of the multinomial distribution. The following three
components are used to describe it:
\begin{itemize}
\item Distributional component: A random sample of size $n$, $\bold{y}_1,\ldots,\bold{y}_n$, is selected from a
multinomial distribution with parameters $(\piv_i,n_i);\,\piv_i = (\pi_{i1},\ldots,\pi_{iq}),\, i = 1,\ldots,n$.
The density function of $\bar{\textbf{y}}_i={\textbf{y}}_i/n_i$ also
called the scaled multinomial distribution \citep[p~76]{fahrmeirtutz_2001} is,
\begin{equation}s(\bar{\textbf{y}}_i|\boldsymbol\theta_i,\phi,\omega)=\exp\left\{\frac{[\bar{\textbf{y}}'_i\boldsymbol
\theta_i-b(\boldsymbol\theta_i)]}{\phi}\omega_i+c(\textbf{y}_i,\phi,\omega_i)\right\},\label{smd}
\end{equation}
where 
$\boldsymbol\theta_i=\left[\log(\frac{\pi_{i1}}{1-\sum_{j=1}^q\pi_{ij}}),
\ldots,\log(\frac{\pi_{iq}}{1-\sum_{j=1}^q\pi_{ij}})\right]'$,
$b(\boldsymbol\theta_i)=-\log({1-\sum_{j=1}^q\pi_{ij}})$, $c({\textbf{y}}_i,\phi,\omega_i)=\log\left(
\frac{n_i!}{y_{i1}!\ldots y_{iq}!(n_i-y_{i1}-\ldots-y_{iq})!}\right)$, $\omega_i=n_i$ and $\phi=1$. The total number of observations is $N=\sum_{i=1}^n n_i$.
\item Linear predictor: A $q$ dimensional linear predictor, $\etav(\xv)=\Zv(\xv)\betav$, where
$\Zv(\xv)=\bigoplus_{j=1}^{q}\textbf{f}_j(\textbf{x})$, $\textbf{f}_j(\textbf{x})$ is a known vector function of
$\textbf{x}$, $\boldsymbol\beta=[\boldsymbol\beta_1',\ldots, \boldsymbol\beta_q']'$
 is the $p \times 1$ vector of unknown parameters with the $j$th component, $\boldsymbol\beta_j$, of length $p_j$ and $p=\sum_{j=1}^q p_j$.

\item Parametric link function: 
$\muv=\piv=\hv(\alphav,\etav)$, where $\hv(\alphav,\cdot)=[h_1(\alphav,\cdot),\ldots,h_q(\alphav,\cdot)]'$,
$\alphav_{r\times 1}=[\alphav_1',\ldots,\alphav_q']'$, $\alphav_j$ is of length $r_j$ and $\sum_{j=1}^{q}r_j=r$.

\end{itemize}

\subsection{Proposed form of parametric link function}

Several researchers \citep{1988_stukel,1989_czado} propose the following generalization for
binary response models with a logistic link function
\begin{equation*}\mu(\xv)=E(y|\xv)=h(\alphav,\eta)=\frac{\exp\{G(\alphav,\eta)\}}{[1+\exp\{G(\alphav,\eta)\}]},\end{equation*}
where $G(\alphav,\cdot)$ is a generating family with the
unknown link parameter $\alphav$. For example using the generating family by \cite{1989_czado} we get
\begin{equation*}G(\alphav,\eta)
=\left\{\begin{matrix}\frac{(1+\eta)^{\alpha_{1}}-1}{\alpha_{1}} & \text{if} & \eta\geq 0\\
-\frac{(1-\eta)^{\alpha_{2}}-1}{\alpha_{2}} & \text{if} & \eta<0,
\end{matrix}\right.\label{czado}
\end{equation*}
where $\alphav=[\alpha_1,\alpha_2]'$.
Usually when modeling the mean in a multinomial response model the multivariate
 version of the logit model \citep[pp.~267-274]{agresti_2002} is used,
\begin{eqnarray*}\pi_{j}=h_j(\etav)=\frac{\exp(\eta_{j})}{1+\sum_{l=1}^q\exp(\eta_{l})},
\text{ for $j=1,\ldots,q$}.\end{eqnarray*}
An alternative form of the above model is given by using the link function $\gv$ where $\gv=\hv^{-1}$,
\begin{eqnarray*}\eta_j=g_j(\muv)=\log\frac{\mu_{j}}{1-\sum_{j=1}^q\mu_{j}}.\end{eqnarray*}
Analogous to the binary case we propose the following generalization of the multicategorical logit model,
\begin{eqnarray}\pi_{j}=h_j(\alphav,\etav)=\frac{\exp\{G_j(\alphav_j,\eta_j)\}}
{1+\sum_{l=1}^q\exp\{G_l(\alphav_l,\eta_l)\}},\text{ for $j=1,\ldots,q$},\label{Gh}\end{eqnarray} where
$\Gv=[G_1,G_2,\ldots,G_q]'$, $G_j(\alphav_j,\cdot)$ is a generating family for binary response models as
described above,
$h_j(\alphav,\cdot)$ is the $j$th component of $\hv(\alphav,\cdot)=\gv^{-1}(\alphav,\cdot)$. The family
$\Lambda=\{\hv(\alphav,\cdot):\alphav\in\Omegav\}$ includes the multivariate logistic link function if
there exists an
$\alphav_0\in\Omegav$ such that $\Gv(\alphav_0,\cdot)$ is a identity function.


\subsection{Parameter estimation}

Summing up the previous sections we can write the multinomial model with
a parametric link function as,
\begin{equation}\piv(\xv)=\hv[\alphav,\etav(\xv)],\end{equation}
where $\etav(\xv)=[\etav_1(\xv),\ldots,\etav_q(\xv)]' =\Zv(\xv)\betav,\,\xv\in R^{k}$,
 $\boldsymbol\beta$ is a $p \times 1$ vector of unknown parameters
  and $\alphav$ is $r\times 1$ vector of unknown link parameters.
Also,
\begin{equation}\etav_{j}=\etav_j(\xv)=\fv_j(\xv)\betav_j=\gv_j[\alphav,\piv(\xv)],\text{ for $j=1,\ldots,q$},\end{equation}
  where $\gv=[\gv_1,\ldots,\gv_q]'$ is the inverse of $\hv=[\hv_1,\ldots,\hv_q]'$.
  We use the notation $\deltav$ to denote the joint vector of the regression and link parameters, thus $\deltav=[\betav',\alphav']'$ is a vector of length $(p+r)$. The parameter vector $\deltav$
  is estimated using the maximum likelihood estimation (MLE) method.
  A brief description of the procedure is given as follows:

Using the scaled version of the multinomial distribution as described
in equation (\ref{smd}) the log-likelihood function for the sample
$\yv_1,\ldots,\yv_n$ is,
  \begin{eqnarray}l(\deltav) &=& \sum_{i=1}^n l_i(\deltav)\nonumber\\
  &=& \sum_{i=1}^n[\yvb_i'\thetav_i-b(\thetav_i)]n_i+constant.\label{ll}
  \end{eqnarray}
Thus the score function is \citep[p~436]{fahrmeirtutz_2001},
  \begin{eqnarray}\frac{\partial l(\deltav)}{\partial\deltav} &=& \frac{\partial}{\partial\deltav}\sum_{i=1}^n[\yvb_i'\thetav_i-b(\thetav_i)]n_i\nonumber\\
  &=& \sum_{i=1}^n\frac{\partial\muv_i}{\partial\deltav}[Var(\yvb_i)]^{-1}(\yvb_i-\muv_i),
  \label{dlddeltav}
  \end{eqnarray}
  and \citep[p~436]{fahrmeirtutz_2001}
  \begin{eqnarray}-\frac{\partial^2 l(\deltav)}{\partial\deltav\partial\deltav'} &=& \sum_{i=1}^n\frac{\partial\muv_i}{\partial\deltav}[Var(\yvb_i)]^{-1}\frac{\partial\muv_i}{\partial\deltav'}-
  \sum_{i=1}^n\sum_{j=1}^q\frac{\partial^2\theta_{ij}}{\partial\deltav\partial\deltav'}(\bar{y}_{ij}-\mu_{ij})n_i\nonumber\\
  &=& \Hv_n,\ (say).\label{dl2ddeltav}
  \end{eqnarray}
  From equation (\ref{dl2ddeltav}) we get the Fisher information matrix to be
  \begin{eqnarray}\Jv_n &=& -E\left[\frac{\partial^2 l(\deltav)}{\partial\deltav\partial\deltav'}\right] 
  = \sum_{i=1}^n\frac{\partial\muv_i}{\partial\deltav}[Var(\yvb_i)]^{-1}\frac{\partial\muv_i}{\partial\deltav'}.\label{jv}
  \end{eqnarray}

  For maximizing the log-likelihood the Fisher scoring iteration method is used which yields,
\begin{equation} \deltav^{(m+1)}=\deltav^{(m)}+\Jv_n^{-1}\frac{\partial l[\deltav^{(m)}]}{\partial\deltav},\end{equation}
 where $m$ indicates the $m$th iteration. It is also possible to use the
 method given in \cite{1988_stukel} for finding an approximate MLE of $\deltav$. In this method,
 the MLE of $\betav$ is obtained
 by fixing $\alphav$ and denoted by $\hat{\beta}(\alphav)$. An approximate MLE
 of $\deltav$ is then given by $\hat{\deltav}=[\hat{\betav}'(\hat{\alphav}),\hat{\alphav}']'$ which maximizes
 the log-likelihood function over a set $\alphav$.

 \subsection{Asymptotic distribution of $\hat{\deltav}$}\label{asym}

Suppose $\deltavh=[\betavh',\alphavh']'$ denotes the MLE of  $\deltav$. Using equation (\ref{dlddeltav}) and the
central limit theorem we know
$\frac{\partial l(\deltav)}{\partial\deltav}$
   asymptotically follows a
  normal distribution with mean $\mathbf{0}$ and
  variance $\Jv_n$.  By first order Taylor series expansion,
\begin{eqnarray*}\mathbf{0}=\frac{\partial l(\deltavh)}{\partial\deltav}&=&\frac{\partial l(\deltav)}{\partial\deltav}+
  \left[\frac{\partial^2 l(\deltav)}{\partial\deltav\partial\deltav'}\right](\deltavh-\deltav).
  \end{eqnarray*}
  This implies \citep[p~439]{fahrmeirtutz_2001},
  \begin{eqnarray*}\sqrt{N}(\deltavh-\deltav)&=&
  \sqrt{N}\Hv^{-1}_n\frac{\partial l(\deltav)}{\partial\deltav}\nonumber=\sqrt{N}\Jv_n^{-1}\frac{\partial l(\deltav)}{\partial\deltav}+O_p(N^{-1/2}).\label{tsd}
  \end{eqnarray*}
  Thus, we get that $\deltavh$
  has an asymptotic normal distribution with mean
  $\deltav$ and variance $\Jv_n^{-1}$.

\section{Orthogonalization of link and regression parameter vectors}\label{ortho}

In this section we discuss certain conditions for which the link parameters are approximately orthogonal to the
regression parameters in a neighbourhood asymptotically. In our numerical examples we show that approximate
orthogonality of the parameters reduces the variance inflation of $\bvh$  while increasing the numerical stability of the computations.
The family of link functions
for which the regression parameters are approximately orthogonal to the link parameters in a neighbourhood are also  location and scale invariant. \citet{1989_li} noted the importance of a family of link function being location and scale invariant.
In their paper they observed that for a unspecified link function
the intercept parameter was not identified while the slope parameter was identified only up to a multiplicative
constant. Thus any variation in the location and scale was absorbed by the link function.
\\

\textbf{Proposition 1:} The regression parameter vector $\betav$ and link parameter vector $\alphav$ are
approximately orthogonal in a neighbourhood around $\etav_0$, if the family of link functions $\hv(\alphav,\cdot)$
satisfies the following conditions, (i) there exists $\etav_0$ and $\piv_0$ such that
\begin{equation}\hv(\alphav,\etav_0)=\piv_0,\ \forall\ \alphav\in\Omega, \label{cond1}\end{equation} and (ii)
there exists
a $\sv_0$ such that
\begin{equation}\frac{\partial \hv(\alphav,\etav)}{\partial\etav}\left|_{\etav=\etav_0}\right.=\sv_0,\
\forall\ \alphav\in\Omega.\label{cond2}\end{equation}

\textbf{Proof:} By first order Taylor series expansion of $\hv(\alphav,\etav)$ around $\etav_0$ and equations
(\ref{cond1}) and (\ref{cond2}),
\begin{eqnarray}\hv(\alphav,\etav) &\approx& \hv(\alphav,\etav_0)+\frac{\partial \hv(\alphav,\etav)}{\partial\etav}
\left|_{\etav=\etav_0}\right.(\etav-\etav_0)\label{te}\\&=&
\piv_0+\sv_0(\etav-\etav_0).\label{ute}\end{eqnarray} Equation (\ref{ute}) shows that the family of link
functions $\hv(\alphav,\cdot)$ is approximately independent of $\alphav$ in a neighbourhood of $\etav_0$
where approximation (\ref{te}) holds asymptotically. Hence, if the conditions (\ref{cond1}) and (\ref{cond2}) are
satisfied, then the regression parameters $\betav$ and the link parameters $\alphav$ are
approximately orthogonal in a neighbourhood of $\etav_0$ asymptotically.\\\par
Extending the definition of a location and scale
invariant family given by \cite{1997_czado} to the multiple response case we state: a family $\Lambdav$ is said
to be {\it location and scale invariant} if for every $\hv\in \Lambdav$, the function
$\hv^*(\alphav,\etav)=\hv(\alphav,\av+\bv\etav)\notin\Lambdav$ for all $\av\neq \mathbf{0}_a$ and $\bv\neq
\mathbf{0}_b$ or $\Iv_b$, where $\av=[a_1,\ldots,a_q]'$, $\bv=diag(b_1,\ldots,b_q)$, $\mathbf{0}_a$ is a matrix
of the same order as $\av$ with all elements zero, $\mathbf{0}_b$ is a matrix of same order as $\bv$ with all
elements zero and $\Iv_b$ is an identity matrix with the same order as $\bv$.\\

\textbf{Proposition 2:} If every member of a family $\Lambdav=\{\hv(\alphav,\cdot):\alphav\in\Omegav\}$ satisfies
conditions (\ref{cond1}) and (\ref{cond2}) for fixed $\piv_0$ and $\sv_0$, then $\Lambdav$ is location and scale invariant.\\

\textbf{Proof: } Suppose every member $\hv$ of the family $\Lambdav$ satisfies conditions (\ref{cond1}) and
(\ref{cond2}) for  fixed $\piv_0$ and $\sv_0$. Define $\hv^*(\alphav,\etav)=\hv(\alphav,\av+\bv\etav)$, then at
$\etav^*=\bv^{-1}(\etav_0-\av)$,
\begin{equation*}\hv^*(\alphav,\etav^*)=\hv(\alphav,\etav_0)=\piv_0\ \forall\ \alphav\ \in\Omegav,\end{equation*}
 where $\av\neq \mathbf{0}_a$ and $\bv\neq
 \mathbf{0}_b$ or $\Iv_b$. Thus, equation (\ref{cond1}) is satisfied by $\hv^*$ at $\etav=\etav^*$.
 Also,
\begin{eqnarray*}\frac{\partial\hv^*(\alphav,\etav)}{\partial\etav}|_{\etav=\etav^*}
 &=& \frac{\partial\hv(\alphav,\av+\bv\etav)}{\partial\etav}|_{\etav=\etav^*}=
\bv\frac{\partial\hv(\alphav,\etav)}{\partial\etav}|_{\etav=\etav_0}= \bv\sv_0\neq\sv_0,
\end{eqnarray*}
for $\bv\neq \Iv_b$, implying $\hv^*\notin\Lambdav$ for $\av\neq \mathbf{0}_a$ and $\bv\neq
 \mathbf{0}_b$ or $\Iv_b$. Hence the family $\Lambdav$ is
location and scale invariant.

\subsection{Construction of $(\piv_0,\sv_0)$-standardized link families at $\etav_0$}

A family $\Lambdav$ satisfying conditions (\ref{cond1}) and (\ref{cond2})
is called $(\piv_0,\sv_0)$-standardized at $\etav_0$ \citep{1997_czado}. \\ 

\textbf{Proposition 3:} Suppose $\Gv(\alphav,\etav)=[G_1(\alphav_1,\eta_1),\ldots,G_q(\alphav_q,\eta_q)]'$ where
$\alphav=[\alphav_1,\ldots,\alphav_q]'$ and $\etav=[\eta_1,\ldots,\eta_q]'$,
such that each \{$G_j(\alphav_j,\eta_j),\ \alphav_j\in\Omega_j$\} is a generating family 
for binary response models and are $(\mu_{0j}, s_{0j})$-standardized
 at $\eta_{0j}$ for $j=1,2,\ldots,q$. Then the family $\Lambdav_{\gv}=\{\Gv(\alphav,\cdot):
 \alphav\in\Omegav=\Omega_1\times\ldots\times\Omega_q$\} is $(\piv_0,\sv_0)$-standardized
 at $\etav_0$, where $\etav_0=[\eta_{01},\ldots,\eta_{0q}]'$, $\piv_0=[\mu_{01},\ldots,\mu_{0q}]'$, and $\sv_0=diag\{s_{01},\ldots,s_{0q}\}$.\\

 \textbf{Proof: } Since, $G_j(\alphav_j,\eta_j)$ is $(\mu_{0j}, s_{0j})$-standardized
 at $\eta_{0j}$,
  $G_j(\alphav_j,\eta_{0j})=\mu_{0j},\ \forall\ \alphav_j\in\Omega_j, $ and
 $\frac{\partial G_j(\alphav,\eta)}{\partial\eta}|_{\eta=\eta_{0j}}= s_{0j},\ \forall\ \alphav_j\in\Omega_j$.
 Thus,
 $\Gv(\alphav,\etav_0) = [\mu_{01},\mu_{02},\ldots,\mu_{0q}]'=\piv_0,\ \forall\ \alphav\in\Omegav$,  and
 $\frac{\partial \Gv(\alphav,\etav)}{\partial\etav}|_{\etav=\etav_0}= diag\{s_{01},s_{02},\ldots,s_{0q}\}=\sv_0,\ \forall\
 \alphav\in\Omegav$.
  Hence, the family $\Lambdav_{\gv}$ is $(\piv_0,\sv_0)$-standardized
 at $\etav_0$.\\

For using a $(\piv_0,\sv_0)$-standardized generating family at $\etav_0$ three parameters, $\piv_0$, $\sv_0$
 and $\etav_0$, need to be estimated.
 Avoiding estimation of extra parameters, the generating family can be standardized by choosing
 $\piv_0=\betav_0$, $\etav_0=\betav_0$ and $\sv_0=\Iv$. This selection allows for a meaningful
interpretation of $\piv_0$ \citep{1997_czado}, when centered covariates are used.
The $(\piv_0=\betav_0,\sv_0=\Iv)$-standardized at $\etav_0=\betav_0=[\beta_{10}\ldots,\beta_{q0}]'$ generating
family  is denoted by $\Gv_c(\alphav,\etav)$, where the $j$th component of $\Gv_c$ is,
\begin{equation*}G_{cj}(\alphav_j,\eta_j)=\beta_{j0}+G(\alphav_j,\eta_{cj}),\ \eta_{cj}=\eta_j-\beta_{j0}.
\end{equation*}
Here, $G(\alphav_j,\cdot)$ is a $(\mu_0=0,\ s_0=1)$-standardized at $\eta_0=0$ generating family
for binary response models and $\beta_{j0}$
is the intercept parameter for the $j$th response.

Using condition (\ref{Gh}), the family of link functions, $\hv(\alphav,\etav)$, for the multinomial response
model is $(\piv_0,\sv_0)$-standardized at $\etav_0$ when
$\pi_{0j}=\frac{\exp(\beta_{j0})}{1+\sum_{j=1}^q\exp(\beta_{j0})}$, and $\sv_0=\cv\Iv$, $\cv$ is
 a constant matrix with its $(j,k)$th element equal to,
     \begin{equation*} c_{jk}=\left\{\begin{matrix} \pi_{0j}(1-\pi_{0j}) & \text{if} & j=k\\ -\pi_{0j}\pi_{0k} & \text{if} & j\neq k.\end{matrix}\right.\end{equation*}
As an example if we are using the generating family as suggested by \cite{1989_czado}
for binary response models as our $G_j$, the
generating family for multinomial responses  is then $(\piv_0=\betav_0,\sv_0=\Iv)$-standardized at
$\etav_0=\betav_0$
\begin{equation}G_{cj}(\alphav_j,\eta_j)=\betav_{j0}+\left\{\begin{matrix}\frac{(1+\eta_{cj})^{\alpha_{j1}}-1}{\alpha_{j1}} & \text{if} & \eta_{cj}\geq 0\\
-\frac{(1-\eta_{cj})^{\alpha_{j2}}-1}{\alpha_{j2}} & \text{if} & \eta_{cj}<0,
\end{matrix}\right.\label{etabeta0}
\end{equation}
 where $\eta_{cj}=\eta_{j}-\beta_{j0}$, for $j=1,\ldots,q$.

In our numerical example (given in Section \ref{example}) we observe that
  the variance inflation ratios are reduced when a $(\piv_0=\betav_0,\sv_0=\Iv)$-standardized
  generating family at
$\etav_0=\betav_0$ is used instead of using a generating family which is
$(\piv_0=\mathbf{0},\sv_0=\Iv)$-standardized at $\etav_0=\mathbf{0}$. Also, we observe that the Newton-Raphson
iteration method does not converge when the $(\piv_0=\mathbf{0},\sv_0=\Iv)$-standardized generating family at
$\etav_0=\mathbf{0}$ is selected, and the grid selection method of \cite{1988_stukel} has to be implemented. For
the generating family $(\piv_0=\betav_0,\sv_0=\Iv)$-standardized at $\etav_0=\betav_0$ the Newton-Raphson
algorithm however converges. The estimation of unknown parameters requires less computational time when the
Newton-Raphson method converges instead of using the grid searching method. Hence, we are able to show in our
example that the numerical stability is increased and the computational time is reduced when using a
$(\piv_0=\betav_0,\sv_0=\Iv)$-standardized at $\etav_0=\betav_0$ generating family .

 \section{Interval estimation of the percentiles}\label{cr}
Suppose we define $\Sv_{\piv_0}(\deltav)$ as the settings
of the control variables at which $\piv_0=\hv[\alphav,\etav(\xv)]$,
\begin{equation}\Sv_{\piv_0}(\deltav)=\{\xv\in R^k:\piv_0=\hv[\alphav,\etav(\xv)]\}.\end{equation}
Then, $\Sv_{\piv_0}(\deltav)$ can be called the
$\piv_0$th percentile of the multinomial distribution. In this section we propose three different methods for determining confidence regions for
$\Sv_{\piv_0}(\deltav)$.

\subsection{Method 1: an asymptotic conservative confidence region based on ML estimates}\label{crm1}

From Section \ref{asym} we know that $\sqrt{N}(\deltavh-\deltav)\sim MVN(\mathbf{0},\Sigmav(\deltav)=N\Jv_n^{-1})$ asymptotically. 
 This implies, $\sqrt{N}(\betavh_j-\betav_{j})$ follows an asymptotic
 multivariate normal distribution with mean $\mathbf{0}_{p_j}$
and variance $\Sigmav_j$, here $\Sigmav_j$ is a sub matrix
  of $\Sigmav(\deltav)$ corresponding to $\betavh_j$. Using the
  normality of $\sqrt{N}(\betavh_j-\betav_j)$ we get,
  $N(\betavh_j-\betav_{j})'\Sigmav^{-1}_j(\betavh_j-\betav_{j})\sim\chi^{2}_{p_j}$ asymptotically.
   Thus, $Pr[N(\betavh_j-\betav_{j})'\Sigmav^{-1}_j(\betavh_j-\betav_{j})\leq\chi^2_{p_j,(1-\tau)}]=(1-\tau)$,
  where $\chi^2_{p_j,(1-\tau)}$ is the $(1-\tau)$th quantile of a $\chi^{2}_{p_j}$ distribution.
  Using Cauchy Schwarz inequality,
  \begin{eqnarray}\underset{\xv\in R^{k}}{sup}\frac{N[\fv_j(\xv)(\betavh_j-\betav_{j})]^2}{\fv_j(\xv)\Sigmav_j\fv_j'(\xv)}
  &\leq& \underset{\zv\in R^{p_j}}{sup}\frac{N[\zv'(\betavh_j-\betav_{j})]^2}{\zv'\Sigmav_j \zv}\nonumber\\
  &=& N(\betavh_j-\betav_{j})'\Sigmav^{-1}_j(\betavh_j-\betav_{j}),\nonumber\\
  &&
  \end{eqnarray}
  thus,
  \begin{equation}\frac{N[\fv_j(\xv)(\betavh_j-\betav_{j})]^2}{\fv_j(\xv)\Sigmav_j\fv_j'(\xv)}\leq N(\betavh_j-\betav_{j})'\Sigmav^{-1}_j(\betavh_j-\betav_{j})
  \text{ for all $\xv\in R^k$}.\label{ci}\end{equation}
  Suppose we define two events $A_j$ and $B_j$ as,
  $A_j=\left[\frac{N[\fv_j(\xv)(\betavh_j-\betav_{j})]^2}
  {\fv_j(\xv)\Sigmav_j\fv_j'(\xv)}\leq\right.$
   $\left.\chi^2_{p_j,(1-\tau)},\ \forall\ \xv\in R^k\right]$ and $B_j=\left[N(\betavh_j-\betav_{j})'\Sigmav^{-1}_j(\betavh_j-\betav_{j})\leq\chi^2_{p_j,(1-\tau)}\right]$.
From equation (\ref{ci}), we know that $B_j\subset A_j$, thus,
\begin{equation*}P[\eta_j(\xv)\in \Cv_j(\xv),\ \forall\ \xv \in R^k]\geq (1-\tau),\text{ for all $j=1,2,\ldots,q$},
  \end{equation*}
  where $\mathcal{C}_j(\xv)=\{\xi\in R:L_j(\xv)\leq \xi\leq U_j(\xv)\}$,
\begin{eqnarray}L_j(\xv) &=& \fv_j(\xv)\betavh_j
  -\sqrt{N^{-1}\fv_j(\xv)\Sigmav_j\fv'_j(\xv)\chi^2_{p_j,(1-\tau)}},\nonumber\\
  U_j(\xv)&=&\fv_j(\xv)\betavh_j
  +\sqrt{N^{-1}\fv_j(\xv)\Sigmav_j\fv'_j(\xv)\chi^2_{p_j,(1-\tau)}}.\label{ljuj}
  \end{eqnarray}
  Then, using Boole's inequality,
\begin{eqnarray*}&&Pr[\eta_j(\xv)\in \Cv_j(\xv),\ \forall\ \xv \in R^k,\text{ for all $j=1,2,\ldots,q$}]\geq (1-q\tau),
\end{eqnarray*}
which implies,
\begin{eqnarray} Pr[\etav(\xv)\in \Cv(\xv),\ \forall\ \xv \in R^k]\geq (1-q\tau),\label{cifetav}
  \end{eqnarray}
  where $\Cv(\xv)=\times_{j=1}^q\Cv_j(\xv)$.
If we now denote $P_{L,j}(\xv)=\min_{\xiv\in\mathcal{C}(\xv)}h_j(\alphav,\xiv)$ and
$P_{U,j}(\xv)=\max_{\xiv\in\mathcal{C}(\xv)}h_j(\alphav,\xiv)$, for $j=1,2,\ldots,q$,
then using the result given in \citep[p~240]{1973_rao},
  \begin{eqnarray*}&&Pr[P_{L,j}(\xv)\leq h_j(\alphav,\etav)\leq P_{U,j}(\xv), \forall\ \xv\in R^k \text{ and }\forall\
j=1,\ldots,q\nonumber] \geq (1-q\tau),
\end{eqnarray*}
This implies that,
\begin{eqnarray}Pr[\Pv_L(\xv)\leq \hv\{\alphav,\etav(\xv)\}\leq
\Pv_U(\xv), \forall\ \xv\in R^k]
\geq (1-q\tau),\label{cifpiv}
  \end{eqnarray}
where $\Pv_L$ and $\Pv_U$ are $q$ dimensional vectors with their $j$th elements equal to $P_{L,j}$ and
$P_{U,j}$, respectively. Since the link parameter $\alphav$ is unknown, we use $\alphavh$ (MLE of $\alphav$)
for computing $P_{L,j}$ and
$P_{U,j}$, for $j=1,2,\ldots,q$.

 Estimating ${\Sv_{\piv_0}}(\deltav)$ by ${\Sv_{\piv_0}}(\hat{\deltav})$, we get a approximate $100(1-\tau')\%$
 ($\tau'=q\tau$) conservative confidence region for $\Sv_{\piv_0}(\deltav)$ as
\begin{eqnarray}\{\zetav\in R^k: \Pv_L(\xv) &\leq& \hv[\alphavh,\etah(\zetav)]
\leq \Pv_U(\xv)\text{ for all $\xv\in{\Sv_{\piv_0}}(\hat{\deltav})$}\}.\label{crco}
\end{eqnarray}

\subsection{Method 2: confidence region using the likelihood ratio test}\label{lrt}
We derive the confidence region of $\Sv_{\piv_0}(\deltav)$ using the likelihood ratio (LR) test corresponding to the
hypotheses, $H_0: \xv\in\Sv_{\piv_0}(\deltav)$ versus $H_1:\xv\notin\Sv_{\piv_0}(\deltav)$.
Under the null hypothesis we have,
\begin{eqnarray*}\etav(\xv) &=& \gv(\alphav,\piv_0),\end{eqnarray*}
which implies
\begin{eqnarray*} \fv_j(\xv)\betav_{j}&=&\gv_j(\alphav,\piv_0)
\ for\ j=1,2,\ldots,q.
\label{lreq1}
\end{eqnarray*}
Suppose $D(\xv)$ is the deviance \citep[p~108]{fahrmeirtutz_2001} under the null hypothesis while $D(\xvh)$
 is the deviance
 of
the fitted model. Then, the LR statistic $L(\xv)=D(\xv)-D(\xvh)$ has an asymptotic $\chi^2$ distribution with
$q$ degrees of freedom. The $100(1-\tau)\%$ confidence region for $\Sv_{\piv_0}(\deltav)$ using the LR statistic is
given by
\begin{equation}\{\xv\in R^k: L(\xv)\leq \chi^2_{q,1-\tau}\}.\label{crlr}\end{equation}

\par
\subsection{Method 3: confidence region using the score test }\label{crsc} Suppose,
$\betav_{0}=[\beta_{10},\ldots,\beta_{q0}]'$ and
$\uv_0=\left[\frac{\partial l}{\partial\betav_0}\right]_{\deltavh_0}$. Let $\Sigmavh_0$ be the estimated variance
of $\uv_0$ at $\deltav=\deltavh_0$, where $l(\deltav)$ is the log-likelihood function and $\deltavh_0$ is the MLE of
$\deltav$ under $H_0$ in Section \ref{lrt}. Then $s(\xv)=\uv'_0\Sigmavh_0^{-1}\uv_0$ has an asymptotic $\chi^2$
distribution with $q$ degrees of freedom \citep[p~48]{fahrmeirtutz_2001}.

Using the score test, the $100(1-\tau)\%$ confidence region for $\Sv_{\piv_0}(\deltav)$ is given by
\begin{equation}\{\xv\in R^k: s(\xv)\leq \chi^2_{q,1-\tau}\}.\label{crs}\end{equation}

\section{Example}\label{example}
 We consider a data set based on a combination drug experiment reported by \citep[pp.~429-451]{gennings_1994}. The main goal of the experiment is to study and model the relationship between the dose levels of two
drugs, morphine sulfate and $\Delta^9$-tetrahydro-cannabinol ($\Delta^9$-THC), on the pain relief and toxic
responses of male mice. Eighteen groups of male mice (six animals per group) were randomly assigned to receive
the treatment combinations and three responses were recorded. They were, $E$: the number of mice in each group
who exhibit only pain relief (no toxic effect), $T$: the number of mice in each group experiencing toxic effects
(irrespective of pain relief), and $W$: the number of mice who experienced neither pain relief nor any toxic
effects. So we may consider $E$ and $T$ as the efficacy and toxicity responses of the two analgesic drugs. The
dose levels of the two drugs formed a $3\times 6$ factorial design, where a treatment combination consisted of a
single injection using one of three levels of morphine sulfate (2, 4, 6 mg/kg) in addition to one of 6 levels of
$\Delta^9$-THC (0.5, 1.0, 2.5, 5.0, 10.0, 15.0 mg/kg). The centered dose levels of the two drugs morphine sulfate
and $\Delta^9$-THC were denoted as $x_1$ and $x_2$. The $3\times 6$ factorial design $D$ and the three responses
are given in Table \ref{responses}. Since the responses were obtained from the same mouse they may be correlated.
The binary nature of the responses allowed us to model them using a two category multinomial model. The response
$W$ was taken correlated binary to be the dummy category. For more details on modeling
 binary responses using the multinomial distribution see \citep{mukhopadhyaykhuri_2008b}.

\begin{table}[h]
 \centering{\small
\caption{Design $D$ and responses $\yv_i=[E_i, T_i, W_i]'$. There are $n_i=6$
experimental units for each run.}
\begin{tabular}{cccccc}
\hline\\
\multicolumn{2}{c}{$D$} & {} & \multicolumn{3}{c}{Responses} \\
\cline{1-2}\cline{4-6}\\
 $x_1$ & $x_2$  & {} & $E_i$ & $T_i$ & $W_i$ \\
 \hline\\
  -1.0 & -0.713 & {} & 5 & 0 & 1\\
  -1.0 & -0.646 & {} & 2 & 0 & 4\\
   -1.0 & -0.433 & {} & 2 & 0 & 4\\
  -1.0 & -0.093 & {} & 4 & 1 & 1\\
  -1.0 & 0.597 & {} & 5 & 1 & 0\\
  -1.0 & 1.287 & {} & 3 & 3 & 0\\
  0.0 & -0.713 & {} & 5 & 0 & 1\\
  0.0 & -0.643 & {} & 6 & 0 & 0\\
  0.0 & -0.433 & {} & 5 & 1 & 0\\
  0.0 & -0.093 & {} & 3 & 3 & 0\\
  0.0 & 0.597 & {} & 3 & 3 & 0\\
  0.0 & 1.287 & {} & 3 & 3 & 0\\
  1.0 & -0.713 & {} & 6 & 0 & 0\\
  1.0 & -0.643 & {} & 6 & 0 & 0\\
  1.0 & -0.433 & {} & 6 & 0 & 0\\
  1.0 & -0.093 & {} & 6 & 0 & 0\\
  1.0 & 0.597 & {} & 1 & 5 & 0\\
  1.0 & 1.287 & {} & 0 & 6 & 0\\

 \hline

\end{tabular}\label{responses}}
\end{table}
\subsection{Fitting a generalized multinomial model}
We start by fitting a multinomial regression model with the multicategorical logit link function to the data. The
model is given by
\begin{eqnarray}\eta_{E}(\textbf{x}) &=& \beta_{10}+\beta_{11}x_1+\beta_{12}x_2,\nonumber\\
\eta_{T}(\textbf{x}) &=& \beta_{20}+\beta_{21}x_1+\beta_{22}x_2.\label{aeta}
\end{eqnarray}
The maximum likelihood estimates (MLEs) of $\betav$ and their standard errors are reported in Table
\ref{mleofbeta}. The scaled deviance for the above fitted model is $29.6048$ with 12 degrees of freedom
(p-value=0.0032). Since the p-value is 0.0032, the results show evidence of lack of fit.
We thus consider the proposed parametric link functions for the multinomial
models and see if it is possible to improve the fit. We use two choices for $\etav_0$, the fixed choice of $\etav_0=\mathbf{0}$ and later $\etav_0=\betav_0$.
\begin{table}[h]
 \centering{\small
\caption{Maximum likelihood estimates and standard errors of the parameters in model (\ref{aeta}).}
\begin{tabular}{cccc}
\cline{1-4}
 \multicolumn{1}{c}{Parameter} & {Estimate}  & {Std. error} & {p-value}\\
 \cline{1-4}
 $\beta_{10}$ & 4.5275 & 1.2203 & 0.0002\\
$\beta_{11}$ & 2.9644 & 1.0566 & 0.0050\\
$\beta_{12}$ & 2.5160 & 1.2987 & 0.0527\\
$\beta_{20}$ & 2.8197 & 1.2582 & 0.0250\\
$\beta_{21}$ & 3.6704 & 1.1231 & 0.0011\\
$\beta_{22}$ & 4.7535 & 1.3761 & 0.0006\\

\cline{1-4}
Deviance=29.6048. & & &\\
\cline{1-4}

\end{tabular}\label{mleofbeta}}
\end{table}

The multinomial model with a parametric link function considering $\etav_0$ to be fixed at $\mathbf{0}$ is given
by
\begin{eqnarray}\pi_{ij}=h_j(\alphav,\etav_i)=\frac{\exp\{G_j(\alphav_j,\eta_{ij})\}}
{1+\sum_{l=1}^2\exp\{G_l(\alphav_l,\eta_{il})\}},\text{ for $j=E,T$},\label{Ghfe}\end{eqnarray}
where \citep{1989_czado}
\begin{equation}G_j(\alphav_j,\eta_{ij})
=\left\{\begin{matrix}\frac{(1+\eta_{ij})^{\alpha_{j1}}-1}{\alpha_{j1}} & \text{if} & \eta_{ij}\geq 0\\
-\frac{(1-\eta_{ij})^{\alpha_{j2}}-1}{\alpha_{j2}} & \text{if} & \eta_{ij}<0,
\end{matrix}\right.\label{czadogffe}
\end{equation}
where $\alphav_j=[\alpha_{j1},\alpha_{j2}]'$ for $j=E,T$. The above link function becomes equivalent to
the multicategorical logistic link function when $\alphav=[1,1,1,1]'$. Using the score test by
\citep[p~48]{fahrmeirtutz_2001}, we test the hypotheses,
\begin{eqnarray}H_0 &:& \alpha_{jk}=1\ \text{ versus }
H_1 : \alpha_{jk}\neq 1\ \forall\ j=1,2\text{ and }\ k=1,2.\nonumber
\end{eqnarray}
From results of the score tests we observe that the null
hypotheses are rejected for $\alpha_{11}$ and $\alpha_{12}$.
This implies that there is a need to modify both tails of the first response. Stepwise
selection of each link parameter based on the akaike information criterion (AIC) were considered in the score tests.
 For computing the MLE of the parameters, $\deltav=(\betav',\alphav')'$ we use the method detailed in
\cite{1988_stukel} since the Fisher scoring iterative method does not converge for $\etav_0=\mathbf{0}$. The parameter estimates, standard
errors
 and variance inflation ratios are given in Table \ref{mlegl}.
 The computations are done by once considering $\alphav$ fixed in the
 information matrix and later estimating it from the data set.

 \par
 We also consider the parametric link function with $\etav_0=\betav_0$ (refer to equation (\ref{etabeta0})).
Using score tests we again note that the link parameters $\alpha_{11}$ and $\alpha_{12}$ need to be included in
the model. The Fisher scoring iteration method for obtaining MLE of $\deltav$ converges and the results are given
in Table \ref{mlegl}.

From Table \ref{mlegl} we note that the  deviance using parametric link function for a generating family
standardized at $\etav_0=\mathbf{0}$
  is 23.8866 at 10 degrees of freedom, and for a generating family standardized at $\etav_0=\betav_0$
    is 22.9148 at 10 degrees of
  freedom. Since, we have estimated  two extra parameters for using the parametric link function,
 so the difference between the deviances using logistic link function and parametric link function is a $\chi^2$ distribution
 with 2 degrees of freedom \citep[p~49]{fahrmeirtutz_2001}.
 The differences between the deviances using a multivariate logistic link function and
  parametric link function with generating families
 (\ref{czadogffe}) and (\ref{etabeta0})
  are 5.7182 (p-value  0.0573) and 6.69 (p-value 0.0353), respectively.
  This shows that using the
  parametric family of link functions
  with generating family (\ref{etabeta0}), we are able to significantly improve the fit over
  the multicategory logistic link function.
  Also in Table \ref{mlegl}, we report the estimates of parameters and two estimated standard errors
  for each regression parameter for both the generating family. The first one assumes that the link parameters are fixed
  at their estimated values and second one assumes that the link parameters are estimated from the data set. The
  variance inflation ratio is the ratio of the standard error when $\alphav$ is estimated to the standard error when
   $\alphav$ is fixed. From  Table \ref{mlegl} we note that the variance inflation ratios corresponding to
the parametric link function with $\etav_0=\mathbf{0}$ are higher than those corresponding to $\etav_0=\betav_0$.
This implies that we achieve greater numerical stability when using a parametric link function with
$\etav_0=\betav_0$.

\begin{table}[h]
 \centering{\small
\caption{Maximum likelihood estimates, standard errors and variance inflation ratios.}
\begin{tabular}{ccccccc}
\cline{1-7}
 \multicolumn{1}{c}{Parameter} & {} & {$\etav_0=0$}  & {} & {} & {$\etav_0=\betav_0$} & {}\\
 \cline{3-4}\cline{6-7}
 {} & {} & {Estimates} & {Variance} & {} & {Estimates} & {Variance}\\
 {} & {} & (S.E.) & {inflation} & {} & (S.E.) & {inflation}\\

 \cline{1-7}

 \multicolumn{1}{c}{$\beta_{10}$} & {\ } & 17.4866 & {} & {} & 6.3931 & {}\\
 {$\alphav$ fixed} & {} & (6.8881)  & {} & {} & (1.4117) & {}\\
 {$\alphav$ estimated} & {} & (13.6982)  & {1.9887} & {} & (1.6786) & {1.1891}\\

 \multicolumn{1}{c}{$\beta_{11}$} & {\ } & 14.1469 & {} & {} &  17.3919 & {}\\
 {$\alphav$ fixed} & {} & (6.0235)  & {} & {} & ( 8.4866) & {}\\
 {$\alphav$ estimated} & {} & (11.9014)  & {1.9758} & {} & (9.6415) & {1.1361}\\

 \multicolumn{1}{c}{$\beta_{12}$} & {\ } & 21.5053 & {} & {} & 15.0282 & {}\\
 {$\alphav$ fixed} & {} & (9.8690)  & {} & {} & (7.0234) & {}\\
 {$\alphav$ estimated} & {} & (20.2745)  & {2.0544} & {} & (8.0906) & {1.1520}\\

 \multicolumn{1}{c}{$\beta_{20}$} & {\ } & 12.3022 & {} & {} & 5.1825 & {}\\
 {$\alphav$ fixed} & {} & (4.9916)  & {} & {} & (1.5889) & {}\\
 {$\alphav$ estimated} & {} & (16.7199)  & {3.3496} & {} & (1.9165) & {1.2062}\\

 \multicolumn{1}{c}{$\beta_{21}$} & {\ } & 10.8918 & {} & {} & 7.2412 & {}\\
 {$\alphav$ fixed} & {} & (3.9276)  & {} & {} & (2.2177) & {}\\
 {$\alphav$ estimated} & {} & (14.6531)  & {3.7308} & {} & (2.8268) & { 1.2747}\\

 \multicolumn{1}{c}{$\beta_{22}$} & {\ } & 17.7634 & {} & {} & 7.8120 & {}\\
 {$\alphav$ fixed} & {} & (6.5569)  & {} & {} & (1.8852) & {}\\
 {$\alphav$ estimated} & {} & (24.0555)  & {3.6687} & {} & (2.7628) & {1.4656}\\

 \multicolumn{1}{c}{$\alpha_{11}$} & {\ } & 0.9 & {} & {} & 0.57 & {}\\
 {$\alphav$ estimated} & {} & ( 0.2218)  & {} & {} & (0.1827) & {}\\

 \multicolumn{1}{c}{$\alpha_{12}$} & {\ } & -2.9 & {} & {} & 0.35 & {}\\
 {$\alphav$ estimated} & {} & (4.0716)  & {} & {} & (0.1519) & {}\\
 \cline{1-7}
 {} & {} & {Deviance = 23.8866} & {} & {} & {Deviance = 22.9148} & {}\\

\cline{1-7}

\end{tabular}\label{mlegl}}
\end{table}
\subsection{Percentile estimation}

In this section we apply the three methods of interval estimation and find confidence regions for
$\piv_0=[0.75,0.2]'$. Thus, we are interested in jointly estimating the $ED_{75}$ and $LD_{20}$ percentiles
(where ED=Effective Dose and LD=Lethal Dose). For computing confidence intervals, we use MLE of $\deltav$ for generating family $(\piv_0=\betav_0,\sv_0=\Iv)$ standardized at $\etav_0=\betav_0$ as it provides better numerical stability and also smaller variance inflation ratios. The
estimated $\piv_0$th percentile is given by
\begin{eqnarray}\Sv_{\piv_0}(\deltavh) &=& \{\xv\in R^2: \piv_0=\hv(\alphavh,\etah(\xv))\}\nonumber\\
&=& \left\{ [-0.6715, 0.1365]'\right\}= \left\{ \xv_0\right\},\  say,
\end{eqnarray}
which is a singleton set. In our example we have two categories and three regression parameters in each category,
thus $q=2$ and $p_j=3$. For obtaining $95\%$ confidence regions for $\Sv_{\piv_0}(\deltav)$, we choose
$\tau'=q\tau=0.05$, which gives $\tau=0.025$ and $\chi^2_{p_j,(1-\tau)}=9.35,\,j=E,T$.

Since $\Sv_{\piv_0}(\deltavh)=\left\{ \xv_0\right\}$, the approximate $95\%$ conservative confidence region for
$\Sv_{\piv_0}(\deltav)$ using method 1 (Section \ref{crm1}) is given by
\begin{eqnarray}\{\xv=[x_1,x_2]'\in R^2: \Pv_L(\xv_0) &\leq& \hv[\alphavh,\etah(\xv)]
\leq \Pv_U(\xv_0)\}.
\end{eqnarray}
For computing the above region we first need to find the intervals for $\eta_E$ and $\eta_T$. Using equation
(\ref{ljuj}), we have the intervals $[-4.6794,-1.7894]$ and $[1.0720,1.7006]$ for $\eta_E$ and $\eta_T$,
respectively. For calculating the intervals we use $N=108$ (total number of observations), $\betavh_j$ and
$\hat{\Sigmav}_j$ from Tables \ref{mlegl} and \ref{varb} respectively. To get $\Cv(\xv_0)$ we take the cartesian products of the
intervals of $\eta_E$ and $\eta_T$. 
The next step is to compute $\Pv_L(\xv_0)=[P_{L,1}(\xv_0),P_{L,2}(\xv_0)]'$ and
$\Pv_U(\xv_0)=[P_{U,1}(\xv_0),P_{U,2}(\xv_0)]'$ where
$P_{L,j}(\xv_0)=\min_{\xiv\in\mathcal{C}(\xv_0)}h_j(\alphavh,\xiv)$ and
$P_{U,j}(\xv_0)=\max_{\xiv\in\mathcal{C}(\xv_0)}h_j(\alphavh,\xiv)$, for $j=E,T$. For computing the minimum and
maximum of the function $h_j(\alphavh,\xiv)$ over $\Cv(\xv_0)$ we use a MATLAB program called MCS
\citep{1999_Huyer} and we get
$\Pv_L(\xv_0)=[0.6319,0.1182]'$ and $\Pv_U(\xv_0)=[0.8414, 0.3113]'$. 
Hence, 
the $95\%$ conservative confidence region for $\Sv_{\piv_0}(\deltav)$ by method 1 is given by
\begin{eqnarray}\{\xv=[x_1,x_2]'\in R^2: [0.6319,0.1182]'\leq \hv[\alphavh,\etah(\xv)]\leq [0.8414, 0.3113]'\}\label{rbm1}
\end{eqnarray}

\begin{table}[h]
 \centering{\small
\caption{Estimated variance of $\deltavh$, when model (\ref{aeta}) is fitted using a parametric family of link functions
  with generating family (\ref{etabeta0}).}
\begin{tabular}{ccccccccc}
\cline{1-9}
 \multicolumn{1}{c}{} & {$\betah_{10}$} & {$\betah_{11}$}  & {$\betah_{12}$} & {$\betah_{20}$} & {$\betah_{21}$} & {$\betah_{22}$}  & {$\alphah_{11}$} & {$\alphah_{12}$}\\
 \cline{1-9}
 $\betah_{10}$ & 2.8176 & 7.5455 & 6.5334 & 3.0072 & 3.3696 & 2.9374 & 0.1117 & 0.1346\\
$\betah_{11}$ & 7.5455 & 92.9582 & 73.5399 & 9.5858 & 18.6817 & 13.7987 & -0.3194 & -0.6891\\
$\betah_{12}$ & 6.5334 & 73.5399 & 65.4584 & 9.0370 & 16.7079 & 15.5723 & -0.0182 & -0.5365\\
$\betah_{20}$ & 3.0072 & 9.5858 & 9.0370 & 3.6728 & 4.5693 & 4.1738 & 0.1884 & 0.1163\\
$\betah_{21}$ & 3.3696 & 18.6817 & 16.7079 & 4.5693 & 7.9907 & 6.9814 & 0.2904 & 0.0229\\
$\betah_{22}$ & 2.9374 & 13.7987 & 15.5723 & 4.1738 & 6.9814 & 7.6333 & 0.3409 & 0.0393\\
$\alphah_{11}$ & 0.1117 & -0.3194 & -0.0182 & 0.1884 & 0.2904 & 0.3409 & 0.0334 & 0.0138\\
$\alphah_{12}$ & 0.1346 & -0.6891 & -0.5365 & 0.1163 & 0.0229 & 0.0393 & 0.0138 & 0.0231\\

\cline{1-9}


\end{tabular}\label{varb}}
\end{table}
 \begin{figure}
  \includegraphics[width= 9in]{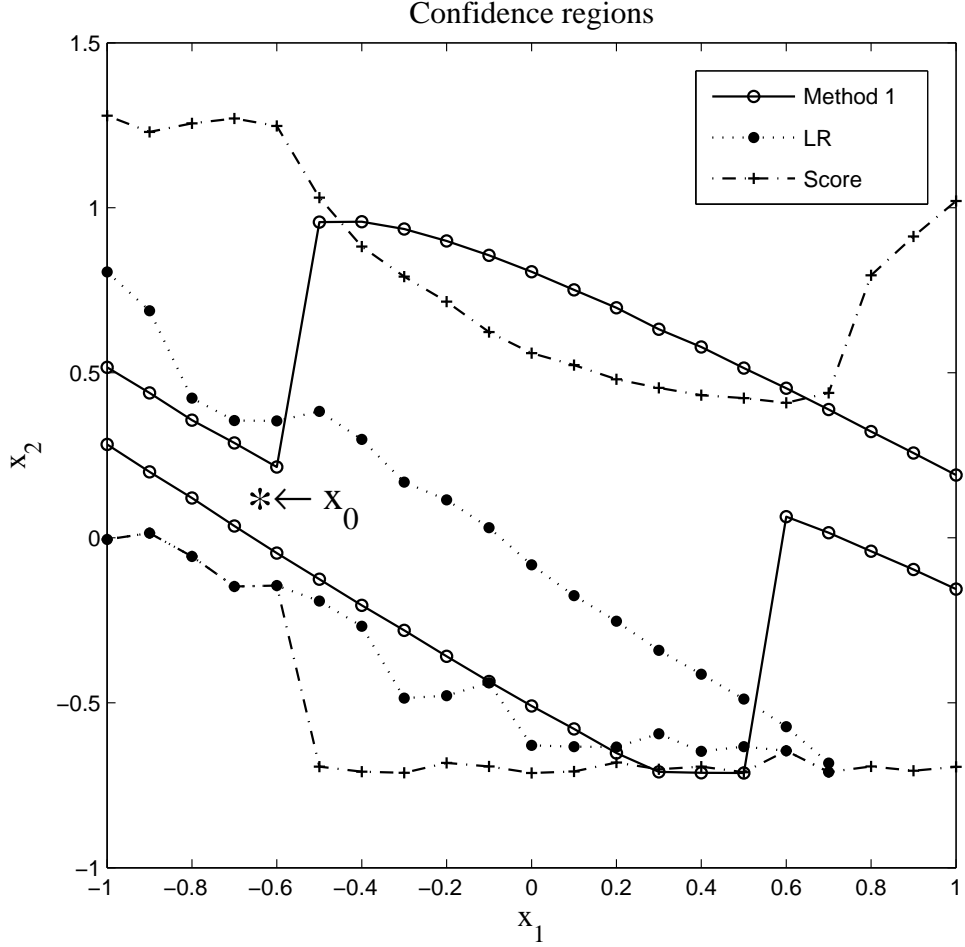}\\
  \caption{Confidence regions for the $\piv_0$th percentile using the three interval estimation methods.}\label{figcr}
\end{figure}
 \par

From equation (\ref{crlr}), the $95\%$ confidence region of $\Sv_{\piv_0}(\deltav)$
 using the LR test (Section \ref{lrt}) is given by
\begin{equation}\{\xv=[x_1,x_2]'\in R^2: L(\xv)\leq  5.99\},\end{equation}
since for $\tau=0.05$, $\chi^2_{2,1-\tau}= 5.99$. The $95\%$ confidence region of $\Sv_{\piv_0}(\deltav)$
using the score test (see equation (\ref{crs})) is,
\begin{equation}\{\xv=[x_1,x_2]'\in R^2: s(\xv)\leq  5.99\},\end{equation}
where $s(\xv)$ is defined as in Section \ref{crsc}.

The confidence regions for the percentiles using the three methods are graphically shown in Figure \ref{figcr}. 
For plotting the confidence regions we choose $21$ values of $x_1$ from the interval $[-1,1]$ at
steps of 0.1.  For each chosen $x_1$, 1000  values of $x_2$ are chosen randomly from  $[-0.7132,1.2868]$. Let
$S_{x_1}$ be the set of the selected $x_2$ values. 
To determine the confidence region by method 1, for each value of $x_1$ we compute
\begin{eqnarray}L_{x_2}(x_1)&=& min\{x_2\in S_{x_1}: [x_1,x_2]' \text{ is in region (\ref{rbm1})} \}\nonumber\\
&& \text{ and }\nonumber\\
U_{x_2}(x_1)&=& max\{x_2\in S_{x_1}: [x_1,x_2]' \text{ is in region (\ref{rbm1})}\}
\end{eqnarray}
Then, by plotting $L_{x_2}(x_1)$ and $U_{x_2}(x_1)$ against $x_1$ we get the lower and upper bounds
 of the region, respectively for method 1. We use the same methodology as above to plot the
 confidence regions by the LR and the score test. 
From Figure \ref{figcr}, we observe that the confidence regions found by using LR test is narrower than the other
two confidence regions, while the score test gives the widest region.

\section{Concluding Remarks} In this article, we have introduced a
family of link functions which is
location and scale invariant and provides local orthogonality between regression
and link parameters for multinomial response models. Using a numerical example
we are able to show that parametric link function provides a better fit
over multivariate logistic link function. We also discussed
three different methods for constructing $100(1-\tau)\%$
confidence regions for the $\piv$th percentile.
\par
 The percentile estimation methods for multinomial models discussed in this article can be used in clinical trials which are involved in determining  dose levels having desired probabilities of both toxicity and efficacy, namely Phase I/II trials \citep{1994_gooley,1998_thall,2001_hughes}. By applying the above interval estimation
 methods experimenters will be able to find confidence regions of dose levels with tolerable toxicity and the desired efficacy.\par
There has been a recent rise of interest among researchers to find designs for logistic regression models which are robust to link misspecification. \citet{2006_biedermann} and \citet{russell_2006} propose robust designs by considering a finite set of plausible link functions while \cite{adewale_2010} uses the family of link functions of
\cite{1981_aranda} in their approach. In the future we plan to use the family of link functions proposed in this article to determine designs for multinomial models which are robust to an incorrectly assumed link function.

\label{conclusion}

\bibliographystyle{elsarticle-harv}
\bibliography{reidview}







\end{document}